\documentclass[journal]{IEEEtran}

\usepackage{cite}
\usepackage{graphicx}
\graphicspath{Figures/}
\usepackage{color}
\DeclareGraphicsExtensions{.pdf}
\usepackage{fixltx2e}
\usepackage{stfloats}
\usepackage{url}
\usepackage{eurosym}
\usepackage[utf8]{inputenc}
\DeclareUnicodeCharacter{20AC}{\euro}
\usepackage{slashbox,booktabs}
\usepackage{multirow}
\usepackage[table]{xcolor}
\usepackage{latexsym}
\usepackage[colorinlistoftodos,prependcaption,textsize=tiny]{todonotes}
\usepackage{mathtools}
\DeclarePairedDelimiter\abs{\lvert}{\rvert}

\usepackage{enumitem}
\setlist[description]{leftmargin=0cm,labelindent=0.6cm}

\begin{document}
%
\title{A Data-driven Bidding Model for a Cluster of Price-responsive Consumers of Electricity}
%
%
%

\author{Javier~Saez-Gallego,
        Juan~M.~Morales,~\IEEEmembership{Member,~IEEE,}
        Marco~Zugno,~\IEEEmembership{Member,~IEEE,}
        Henrik~Madsen 
\thanks{J. Saez-Gallego, J. M. Morales (corresponding author), and H. Madsen are with the Technical University of Denmark, DK-2800 Kgs. Lyngby, Denmark (email addresses: \{jsga, jmmgo, hmad\}@dtu.dk), and their work is partly funded by DSF (Det Strategiske Forskningsr{\aa}d) through the CITIES research center (no. 1035-00027B) and the iPower platform project (no. 10-095378).}
}

%
%

\markboth{ }%
{Saez-Gallego \MakeLowercase{\textit{et al.}}: A Data-driven Bidding Model for a Cluster of Price-responsive Consumers of Electricity}
%



\maketitle

\begin{abstract}
This paper deals with the market-bidding problem of a cluster of price-responsive consumers of electricity. We develop an inverse optimization scheme that, recast as a bilevel programming problem, uses price-consumption data to estimate the complex market bid that best captures the price-response of the cluster. The complex market bid is defined as a series of marginal utility functions plus some constraints on demand, such as maximum pick-up and drop-off rates. The proposed modeling approach also leverages information on exogenous factors that may influence the consumption behavior of the cluster,  e.g.,  weather conditions and calendar effects. We test the proposed methodology for a particular application: forecasting the power consumption of a small aggregation of households that took part in the Olympic Peninsula project. Results show that the price-sensitive consumption of the cluster of flexible loads can be largely captured in the form of a complex market bid, so that this could be ultimately used for the cluster to participate in the wholesale electricity market.
\end{abstract}


%
\IEEEpeerreviewmaketitle

\section{Introduction}


We consider the case of a cluster of flexible power consumers, where \emph{flexibility} is understood as the possibility for each consumer in the cluster to change her consumption depending on the electricity price and on her personal preferences. There are many examples of methods to schedule the consumption of individual price-responsive loads (see, e.g., \cite{Halvgaard2012,ConejoMorales2010,mohsenian2010}).
The portfolio of flexible consumers is managed by a retailer or \emph{aggregator}, which bids in a wholesale electricity market on behalf of her customers. We consider the case where such a market accepts complex bids, consisting of a series of price-energy bidding curves, consumption limits, and maximum pick-up and drop-off rates.
%
In this paper, we present a data-driven methodology for determining the complex bid that best represents the reaction of the pool of flexible consumers to the market price.
The contributions of this paper are fourfold. The first contribution corresponds to the methodology itself: \emph{we propose a novel approach to capture the price-response of a pool of flexible consumers in the form of a market bid using price-consumption data.} In this work, the price is given as the result of a competitive market-clearing process, and we have access to it only from historical records. This is in contrast to some previous works, where the price is treated as a control variable to be decided by the aggregator or retailer. In \cite{Dorini2012}, for example, the relationship between price and consumption is first modeled by a Finite Impulse Response (FIR) function as in \cite{Corradi2012} and peak load reduction is achieved by modifying the price. Similar considerations apply to the works of \cite{Zugno2013,chen2011,meng2013stackelberg, Nikos2015}, where a bilevel representation of the problem is used: the lower-level problem optimizes the household consumption based on the broadcast electricity price, which is determined by the upper-level problem to maximize the aggregator's/retailer's profit.
Another series of studies concentrate on estimating price-energy bids for the participation of different types of flexible loads in the wholesale electricity markets, for example, time-shiftable loads \cite{Mohsenian-Rad}, electric vehicles \cite{Vaya2015} and thermostatically-controlled loads \cite{Li2015}. Contrary to these studies, our approach is data-driven and does not require any assumption about the nature of the price-responsive loads in the aggregation. 

%
%
%
%
	
The second contribution lays in the estimation procedure: \emph{we develop an inverse optimization framework that results in a bilevel optimization problem}.
From a methodological point of view, our approach builds on the inverse optimization scheme introduced in \cite{AhujaInverse}, but with several key differences. First, we let the measured solution  be potentially non-optimal, or even non-feasible, for the targeted optimization problem as in \cite{chan2014generalized} and \cite{Keshavarz_imputinga}. Moreover, we extend the concept of inverse optimization to a problem where the estimated parameters may depend on a set of features and are also allowed to be in the constraints, and not only in the objective function.


Third, \emph{we study heuristic solution methods to reduce the computing times resulting from the consideration of large datasets of features for the model estimation.} We do not solve the resulting bilevel programming problem to optimality but instead we obtain an approximate solution by penalizing the violation of complementarity constraints following a procedure inspired by the work of~\cite{Gabriel2013}. Finally, \emph{we test the proposed methodology using data from a real-world experiment} that was conducted as a part of the \emph{Olympic Peninsula Project} \cite{Olympic}.

\section{Methodology}
%
In this section, we describe the methodology to determine the optimal market bid for a pool of price-responsive consumers. The estimation procedure is cast as a bilevel programming problem, where the upper level seeks to minimize a norm of the estimation error, in particular, the absolute difference between the measured consumption and the estimated one, while the lower level ensures that the estimated consumption is optimal, given the reconstructed bid parameters and the electricity price.


\subsection{Lower-Level Problem: Price-response of the Pool of Consumers}

The lower-level problem models the price-response of the pool of consumers in the form of a market bid, whose parameters are determined by the upper-level problem. The bid is given by $ \theta_{t} = \lbrace a_{b,t}, r^{u}_{t}, r^{d}_{t}, \underline{P}_{t},\overline{P}_{t} \rbrace$, which consists of the marginal utility corresponding to each bid block $b$, the maximum load pick-up and drop-off rates (analogues to the ramp-up and -down limits of a power generating unit), the minimum power consumption, and the maximum power consumption, at time $t \in \mathcal{T} \equiv \lbrace t : t= 1 \ldots T \rbrace$, in that order.
The utility is defined as the benefit that the pool of users obtains from consuming a certain amount of electricity.
The marginal utility $a_{b,t}$ at time $t$ is formed by $b \in \mathcal{B} \equiv \lbrace b: b= 1 \ldots B \rbrace$ blocks, where all blocks have equal size, spanning from the minimum to the maximum allowed consumption. In other words, the size of each block is $\frac{ \overline{P}_t - \underline{P}_t}{B}$. Furthermore, we assume that the marginal utility is monotonically decreasing as consumption increases, i.e., $a_{b,t} \geq a_{b+1,t}$ for all times $t$. Finally, the total consumption at time $t$ is given by the sum of the minimum power demand plus the consumption linked to each bid block, namely, $x_t^{tot} = \underline{P}_{t} + \sum_{b \in \mathcal{B}} x_{b,t}$.

Typically, the parameters of the bid may change across the hours of the day, the days of the week, the month, the season, or any other indicator variables related to the time. Moreover, the bid can potentially depend on some external variables such as temperature, solar radiation, wind speed, etc. Indicator variables and external variables, often referred to as \emph{features}, can be used to explain more accurately the parameters of the market bid that best represents the price-response of the pool of consumers. This approach is potentially useful in practical applications, as numerous sources of data can help better explain the consumers' price-response. We consider the $I$ external variables or features, named $Z_{i,t}$ for $ i \in \mathcal{I} \equiv \lbrace i: i=1, \ldots, I \rbrace$, to be affinely related to the parameters defining the market bid by a coefficient $\alpha_{i}$. This affine dependence can be enforced in the model by letting
$a_{b,t} = a_{b}^{0} + \sum_{i \in \mathcal{I}}  \alpha_{i}^{a}Z_{i,t}$,
$r^{u}_{t} = r^{u0} + \sum_{i \in \mathcal{I}}  \alpha_{i}^{u}Z_{i,t}$,
$r^{d}_{t} = r^{d0} + \sum_{i \in \mathcal{I}}  \alpha_{i}^{d}Z_{i,t}$,
$\overline{P}_{t} = \overline{P}^{0}  + \sum_{i \in \mathcal{I}}  \alpha_{i}^{\overline{P}} Z_{i,t}$,
and $\underline{P}_{t} = \underline{P}^{0} + \sum_{i \in \mathcal{I}} \alpha_{i}^{\underline{P}} Z_{i,t}$. The affine coefficients $\alpha_{i}^{a}$,$\alpha_{i}^{u}$, $\alpha_{i}^{d}$, $\alpha_{i}^{\overline{P}}$ and $\alpha_{i}^{\underline{P}}$, and the intercepts $a_{b}^{0}, r^{u0},r^{d0}, \underline{P}^{0},\overline{P}^{0}$ enter the model of the pool of consumers (the lower-level problem) as parameters, together with the electricity price.

The objective is to maximize consumers' welfare, namely, the difference between the total utility and the total payment:
\begin{subequations} \label{eq:orig}
\begin{align}
&\underset{x_{b,t}}{\rm{Maximize}} \sum_{t \in \mathcal{T}} \left(  \sum_{b \in \mathcal{B}} a_{b,t}x_{b,t} - p_{t}\sum_{b \in \mathcal{B}} x_{b,t} \right)
\end{align}
where $x_{b,t}$ is the consumption assigned to the utility block $b$ during the time $t$, $a_{b,t}$ is the marginal utility obtained by the consumer in block $b$ and time $t$, and $p_t$ is the price of the electricity during time $t$. For notational purposes, let $\mathcal{T}_{-1} = \lbrace t: t=2, \ldots, T\rbrace$. The problem is constrained by
{\allowdisplaybreaks
\begin{IEEEeqnarray}{lLr}
& \underline{P}_{t} + \sum_{b \in \mathcal{B}} x_{b,t} - \underline{P}_{t-1} - \sum_{b \in \mathcal{B}} x_{b,t-1}   \leq r^{u}_{t}
& t \in \mathcal{T}_{-1} \  \   \label{eq:prim1}
\\
& \underline{P}_{t-1} + \sum_{b \in \mathcal{B}} x_{b,t-1}  - \underline{P}_{t}- \sum_{b \in \mathcal{B}} x_{b,t}  \leq r^{d}_{t}
&  t \in \mathcal{T}_{-1} \ \ \label{eq:prim2}
\\
&  x_{b,t} \leq \frac{\overline{P}_{t} - \underline{P}_{t}}{B}
& b \in \mathcal{B}    , t  \in \mathcal{T}   \ \    \label{eq:prim3}
\\
& x_{b,t} \geq 0
& b \in \mathcal{B} ,t  \in \mathcal{T}.  \ \   \label{eq:prim4}
\end{IEEEeqnarray}
}
\end{subequations}
Equations \eqref{eq:prim1} and \eqref{eq:prim2} impose a limit on the load pick-up and drop-off rates, respectively. The set of equations \eqref{eq:prim3} defines the size of each utility block to be equally distributed between the maximum and minimum power consumptions. Constraint \eqref{eq:prim4} enforces the consumption pertaining to each utility block to be positive. Note that, by definition, the marginal utility is decreasing in $x_t$ ($a_{b,t} \geq a_{b+1,t}$), so one can be sure that the first blocks will be filled first. We denote the dual variables associated with each set of primal constraints as  $\lambda_{t}^{u},\lambda_{t}^{d},\overline{\psi}_{b,t}$ and $\underline{\psi}_{b,t}$.

Problem \eqref{eq:orig} is linear, hence it can be equivalently recast as the following set of KKT conditions \cite{luenberger2008linear}, where \eqref{eq:stat1}--\eqref{eq:stat3} are the stationary conditions and \eqref{eq:comp1}--\eqref{eq:comp7} enforce complementarity slackness:
\begin{subequations}\label{eq:KKT}
\begin{IEEEeqnarray}{lLr}
& -\lambda_{2}^{u} + \lambda_{2}^{d} -\underline{\psi}_{b,1} + \overline{\psi}_{b,1}    = a_{b,1} - p_{1}
\qquad \ b \in \mathcal{B} \label{eq:stat1}
\\
&  \lambda_{t}^{u} -\lambda_{t+1}^{u} - \lambda_{t}^{d} + \lambda_{t+1}^{d} -\underline{\psi}_{b,t} + \overline{\psi}_{b,t}   = a_{b,t} - p_{t}
\nonumber
\\
& \qquad\qquad\qquad\qquad\qquad\qquad \ \ \ \ \ \forall b \in \mathcal{B}, t \in \mathcal{T}_{-1} \label{eq:stat2}
\\
&   \lambda_{T}^{u} - \lambda_{T}^{d}   -\underline{\psi}_{b,T} + \overline{\psi}_{b,T} = a_{b,T} - p_{T}
\qquad \ \ b \in \mathcal{B} \label{eq:stat3}
\\
&\underline{P}_{t} + \sum_{b \in \mathcal{B}} x_{b,t} - \underline{P}_{t-1} - \sum_{b \in \mathcal{B}} x_{b,t-1} \leq r^{u}_{t}   \perp \lambda_{t}^{u} \geq 0
 \nonumber
\\
& \qquad\qquad\qquad\qquad\qquad\qquad\qquad\qquad \ \ \   t \in \mathcal{T}_{-1} \label{eq:comp1}
\\
&\underline{P}_{t-1} + \sum_{b \in \mathcal{B}} x_{b,t-1}  -\underline{P}_{t}- \sum_{b \in \mathcal{B}} x_{b,t}\leq r^{d}_{t}   \perp \lambda_{t}^{d} \geq 0
  \nonumber
\\
& \qquad\qquad\qquad\qquad\qquad\qquad\qquad\qquad \ \ \  t \in \mathcal{T}_{-1}  \label{eq:comp2}
\\
& x_{b,t} \leq \frac{ \overline{P}_{t} - \underline{P}_{t}}{B}  \perp \overline{\psi}_{b,t}  \geq 0
\qquad\qquad b \in \mathcal{B},t \in \mathcal{T} \label{eq:comp6}
\\
& 0 \leq x_{b,t}  \perp \underline{\psi}_{b,t} \geq 0
\qquad\qquad\qquad \ \ \  b \in \mathcal{B},t \in \mathcal{T}.\label{eq:comp7}
\end{IEEEeqnarray}
\end{subequations}

\subsection{Upper-Level Problem: Market-Bid Estimation Via Inverse Optimization} \label{sec:upper}
Given a time series of price-consumption pairs $(p_{t},x_{t}^{meas})$, the inverse problem consists in estimating the value of the parameters $\theta_t$ defining the objective function and the constraints of the lower-level problem~\eqref{eq:orig} such that the optimal consumption $x_{t}$ resulting from this problem is as close as possible to the measured consumption $x_{t}^{meas}$ in terms of a certain norm. The parameters of the lower-level problem $\theta_t$ form, in turn, the market bid that best represents the price-response of the pool.

In mathematical terms, the inverse problem can be described as a minimization problem:
\begin{subequations}\label{eq:kktobj_lin1}

\begin{align}
\underset{x,\theta}{\rm{Minimize}} \   \sum_{t \in \mathcal{T}} w_{t}\abs[\Big]{ \underline{P}_{t} + \sum_{b \in \mathcal{B}} x_{b,t} - x_{t}^{meas} } \label{eq:kktobj}
\end{align}
%
subject to
\begin{align}
a_{b,t} &\geq a_{b+1,t} & b \in \mathcal{B}, t \in \mathcal{T} \label{eq:inv31}
\\
&\eqref{eq:KKT} \label{eq:inv41}.
\end{align}
\end{subequations}

Constraints \eqref{eq:inv31} are the upper-level constraints, ensuring that the estimated marginal utility must be monotonically decreasing. Constraints \eqref{eq:inv41} correspond to the KKT conditions of the lower-level problem \eqref{eq:orig}.

Notice that the upper-level variables $\theta_t$, which are parameters in the lower-level problem, are also implicitly constrained by the optimality conditions~\eqref{eq:KKT} of this problem, i.e., by the fact that $x_{b,t}$ must be optimal for~\eqref{eq:orig}. This guarantees, for example, that the minimum power consumption  be positive and equal to or smaller than the maximum power consumption ( $ 0 \leq \underline{P}_{t} \leq \overline{P}_{t}  $). Furthermore, the maximum pick-up rate is naturally constrained to be equal to or greater than the negative maximum drop-off rate  ($-r^{d}_{t} \leq r^{u}_{t} $). Having said that, in practice, we need to ensure that these constraints are fulfilled for all possible realizations of the external variables and not only for the ones observed in the past. We achieved this by enforcing the robust counterparts of these constraints \cite{ben2009robust}. An example is provided in the appendix.


Parameter $w_{t}$ represents the weight of the estimation error at time $t$ in the objective function. These weights have a threefold purpose. Firstly, if the inverse optimization problem is applied to estimate the bid for the day-ahead market, the weights could represent the cost of balancing power at time $t$. In such a case, consumption at hours with a higher balancing cost would weigh more and consequently, would be fit better than that occurring at hours with a lower balancing cost.
Secondly, the weights can include a forgetting factor to give exponentially decaying weights to past observations. Finally, a zero weight can be given to missing or wrongly measured observations.

The absolute value of the residuals can be linearized by adding two extra nonnegative variables, and by replacing the objective equation \eqref{eq:kktobj} with the following linear objective function plus two more constraints, namely, \eqref{eq:inv1} and \eqref{eq:inv2}:
\begin{subequations}\label{eq:kktobj_lin}
\begin{align}
\underset{x_t,\theta_t,e^{+}_{t},e^{-}_{t}}{\rm{Minimize}} \sum_{t=1}^{T} w_{t}( e^{+}_{t} + e^{-}_{t} )  \label{eq:kktobj3}
\end{align}
subject to
\begin{align}
\underline{P}_{t} + \sum_{b \in \mathcal{B}} x_{b,t} - x_{t}^{meas} &= e^{+}_{t} -e^{-}_{t} & t  \in \mathcal{T} \label{eq:inv1}
\\
e^{+}_{t},e^{-}_{t} &\geq 0 &  t \in \mathcal{T} \label{eq:inv2}
\\
a_{b,t} &\geq a_{b+1,t} &  t \in \mathcal{T} \label{eq:inv3}
\\
\eqref{eq:KKT} \label{eq:inv4}.
\end{align}
\end{subequations}

In the optimum, and when $w_t >0$, \eqref{eq:inv1} and \eqref{eq:inv2} imply that $e^{+}_{t}=x_{t} - x_{t}^{meas}$ if $x_{t} \geq x_{t}^{meas}$, else $e^{-}_{t}=x_{t}^{meas} - x_{t}$. By using this reformulation of the absolute value, the weights could also reflect whether the balancing costs are symmetric or skewed. In the latter case, there would be different weights for $e^{+}_{t}$ and $e^{-}_{t}$.

\section{Solution Method} \label{sec:solution}

The estimation problem~\eqref{eq:kktobj_lin} is non-linear due to the complementarity constraints of the KKT conditions of the lower-level problem \eqref{eq:KKT}. There are several ways of dealing with these constraints, for example, by using a non-linear solver \cite{ferris2000complementarity}, by recasting them in the form of disjunctive constraints \cite{fortuny1981representation}, or by using SOS1 variables \cite{beale1970special}. In any case, problem~\eqref{eq:kktobj_lin} is \mbox{NP-hard} to solve and the computational time grows exponentially with the number of complementarity constraints.
Our numerical experiments showed that, for realistic applications involving multiple time periods and/or numerous features, none of these solution methods were able to provide a good solution to problem~\eqref{eq:kktobj_lin} in a reasonable amount of time. To tackle this problem in an effective manner, we propose the following two-step solution strategy, inspired by \cite{Gabriel2013}:




\begin{description}
\item[Step 1:] \  Solve a linear relaxation of the mathematical program with equilibrium constraints~\eqref{eq:kktobj_lin} by penalizing violations of the complementarity constraints.
\item[Step 2:] \ Fix the parameters defining the constraints of the lower-level problem~\eqref{eq:orig}, i.e., $r^{u}, r^{d}, \underline{P},\overline{P}, \alpha_{i}^{a},\alpha_{i}^{d},\alpha_{i}^{\overline{P}}$ and $\alpha_{i}^{\underline{P}}$,  at the values estimated in Step 1. Then, recompute the parameters defining the utility function, $a_{b,t}$ and $\alpha_d^a$. To this end, we make use of the primal-dual reformulation of the price-response model~\eqref{eq:orig} \cite{chan2014generalized}.
\end{description}

Both steps are further described in the subsections below. Note that the obtained solution does not guarantee optimality, and it is only proved to work satisfactorily in practice (see the case study in Section~\ref{sec:studycase}).

\subsection{Penalty Method} \label{sec:penal}


The so-called penalty method is a convex (linear) relaxation of a mathematical programming problem with equilibrium constraints, whereby the complementarity conditions of the lower-level problem, that is, problem~\eqref{eq:orig}, are moved to the objective function \eqref{eq:kktobj3} of the upper-level problem. Thus,
we penalize the sum of the dual variables of the inequality constraints of problem~\eqref{eq:orig} and their slacks, where the slack of a ``$\leq$"-constraint is defined as the difference between its right-hand and left-hand sides, in such a way that the slack is always nonnegative. For example, the slack of the constraint relative to the maximum pick-up rate~\eqref{eq:prim1} is defined as $ s_t = r^{u}_{t} - \underline{P}_{t} - \sum_{b \in \mathcal{B}} x_{b,t} + \underline{P}_{t-1} + \sum_{b \in \mathcal{B}} x_{b,t-1}$, and analogously for the rest of the constraints of the lower-level problem.

The penalization can neither ensure that the complementarity constraints are satisfied, nor that the optimal solution of the inverse problem is achieved.  Instead, with the penalty method, we obtain an approximate solution. In the case study of Section~\ref{sec:studycase}, nonetheless, we show that this solution performs notably well.

After relaxing the complementarity constraints \eqref{eq:comp1}--\eqref{eq:comp7}, the objective function of the estimation problem writes as:

\begin{subequations}\label{eq:penalization}
{\allowdisplaybreaks
\begin{align}
\underset{\mathbf{\Omega}}{\rm{Minimize}} \ & \sum_{t \in \mathcal{T}} w_{t}( e^{+}_{t} + e^{-}_{t})+ \nonumber
\\
&  L  \bigg(
 \sum_{\begin{subarray}{c} b \in \mathcal{B} \\ t \in \mathcal{T} \end{subarray} } w_{t}  \Big( \psi^{\overline{P}}_{b,t} + \psi^{\underline{P}}_{b,t} + \frac{\overline{P}_{t} - \underline{P}_{t}}{B}  \Big)  + \nonumber \\
  &\sum_{t\in \mathcal{T}_{-1}}  w_{t} \Big( \lambda_{t}^{u} + \lambda_{t}^{d} + r^{u}_{t} + r^{d}_{t} \Big) \bigg) \label{eq:objPenalty}
\end{align}
}
{
    \def\OldComma{,}
    \catcode`\,=13
    \def,{%
      \ifmmode%
        \OldComma\discretionary{}{}{}%
      \else%
        \OldComma%
      \fi%
    }%
with the variables being $ \mathbf{\Omega}= \lbrace x_t,\theta_t,e^{+}_{t}, e^{-}_{t} , \psi^{\overline{P}}_{t},\psi^{\underline{P}}_{t}, \lambda_{t}^{u} , \lambda_{t}^{d},
\overline{\phi}_{i,t},  \underline{\phi}_{i,t} , \overline{\varphi}_{i,t}, \underline{\varphi}_{i,t},\overline{\eta}_{i,t}, \underline{\eta}_{i,t} \rbrace $, subject to the following constraints:
}
%
%
%
\begin{IEEEeqnarray}{ll}
& \rm{ \left( \ref{eq:inv1}\right)-\left(\ref{eq:inv3} \right)}
,
 \rm{ \left( \ref{eq:prim1}\right)-\left(\ref{eq:prim4}  \right)}
 ,
 \rm{ \left( \ref{eq:stat1}\right)-\left(\ref{eq:stat3}  \right)}
\\
&\lambda_{t}^{u} , \lambda_{t}^{d} \geq 0 \qquad\qquad t \in \mathcal{T}_{-1}  \label{eq:dualfeas1}
\\
& \overline{\psi}_{b,t},\underline{\psi}_{b,t} \geq 0 \qquad \ \   b \in \mathcal{B},t \in \mathcal{T} \label{eq:dualfeas2}.
\end{IEEEeqnarray}

\end{subequations}

The objective function \eqref{eq:objPenalty} of the relaxed estimation problem is composed of two terms. The first term represents the weighted sum of the absolute values of the deviations of the estimated consumption from the measured one.
The second term, which is multiplied by the penalty term $L$, is the sum of the dual variables of the constraints of the consumers' price-response problem plus their slacks. Note that summing up the slacks of the constraints of the consumers' price-response problem eventually boils down to summing up the right-hand sides of such constraints.
The weights of the estimation errors ($w_t$) also multiply the penalization terms. Thus, the model weights violations of the complementarity constraints in the same way as the estimations errors are weighted.

Objective function \eqref{eq:objPenalty} is subject to the auxiliary constraints modeling the absolute value of estimation errors (\ref{eq:inv1})--(\ref{eq:inv2}); the upper-level-problem constraints imposing monotonically decreasing utility blocks \eqref{eq:inv3}; and the primal and dual feasibility constraints of the lower-level problem, (\ref{eq:prim1})--(\ref{eq:prim4}), (\ref{eq:stat1})--(\ref{eq:stat3}), and (\ref{eq:dualfeas1})--(\ref{eq:dualfeas2}).

The penalty parameter $L$ should be tuned carefully. We use cross-validation to this aim, as described in the case study; we refer to Section \ref{sec:studycase} for further details.

Finding the optimal solution to problem~\eqref{eq:penalization} is computationally cheap, because it is a linear programming problem. On the other hand, the optimal solution to this problem might be significantly different from the one that we are actually looking for, which is the optimal solution to the original estimation problem~\eqref{eq:kktobj_lin}. Furthermore, the solution to~\eqref{eq:penalization} depends on the user-tuned penalization parameter $L$, which is given as an input and needs to be decided beforehand.

\subsection{Refining the Utility Function} \label{sec:refine}

In this subsection, we elaborate on the second step of the strategy we employ to estimate the parameters of the market bid that best captures the price-response of the cluster of loads. Recall that this strategy has been briefly outlined in the introduction of Section~\ref{sec:solution}. The ultimate purpose of this additional step is to re-estimate or refine the parameters characterizing the utility function of the consumers' price-response model~\eqref{eq:orig}, namely, $a_{b}^{0}$ and the coefficients $\alpha_{i}^{a}$. In plain words, we want to improve the estimation of these parameters with respect to the values that are directly obtained from the relaxed estimation problem~\eqref{eq:penalization}.  With this aim in mind, we fix the parameters defining the constraints of the cluster's price-response problem~\eqref{eq:orig} to the values estimated in Step 1, that is, to the values obtained by solving the relaxed estimation problem~\eqref{eq:penalization}. Therefore, the bounds $\underline{P}_t,\overline{P}_t$ and the maximum pick-up and drop-off rates $r^{u}_t, r^{d}_t$
are now treated as given parameters in this step. Consequently, the only upper-level variables that enter the lower-level problem~\eqref{eq:orig}, namely, the intersects $a_{b}^{0}$ of the various blocks defining the utility function and the linear coefficients $\alpha_{i}^{a}$, appear in the objective function of problem~\eqref{eq:orig}. This will allow us to formulate the utility-refining problem as a linear programming problem.

Indeed, consider the primal-dual optimality conditions of the consumers' price-response model~\eqref{eq:orig}, that is, the primal and dual feasibility constraints and the strong duality condition. These conditions are also necessary and sufficient for optimality due to the linear nature of this model.

We determine the (possibly approximate) block-wise representation of the measured consumption at time $t$, $x^{meas}_{t}$, which we denote by $\sum_{b \in \mathcal{B}} x^{m'}_{b,t}$ and is given as a sum of $B$ blocks of size $\frac{ \overline{P}_{t} - \underline{P}_{t}}{B}$ each. In particular, we define $\sum_{b \in \mathcal{B}} x^{meas'}_{b,t}$ as follows:
\begin{align*}
\sum_{b \in \mathcal{B}} x^{m'}_{b,t} =
  \begin{cases}
    \overline{P}_{t}       & \quad \text{if } x^{meas}_t > \overline{P}_{t} \ , \ t \in \mathcal{T}\\
    x^{meas}_t  & \quad \text{if }  \underline{P}_{t} \leq x^{meas}_t \leq \overline{P}_{t} \ , \ t \in \mathcal{T}\\
    \underline{P}_{t}  & \quad \text{if } x^{meas}_t < \underline{P}_{t}\ , \   t \in \mathcal{T}.\\
  \end{cases}
\end{align*}
where each $ x^{m'}_{b,t}$ is determined such that the blocks with higher utility are filled first. Now we replace $x_t$ in the primal-dual reformulation of~\eqref{eq:orig} with $\sum_{b \in \mathcal{B}} x^{m'}_{b,t}$. Consequently, the primal feasibility constraints are ineffective and can be dropped.

Once $x_t$ has been replaced with $\sum_{b \in \mathcal{B}} x^{m'}_{b,t}$  in the primal-dual reformulation of~\eqref{eq:orig} and the primal feasibility constraints have been dropped, we solve an optimization problem (with the utility parameters $a_{b}$ and $\alpha_{i}^{a}$ as decision variables) that aims to minimize the weighted duality gap, as in \cite{chan2014generalized}. For every time period $t$ in the training data set, we obtain a contribution ($\epsilon_t$) to the total duality gap ($\sum_{t \in \mathcal{T}}\epsilon_t$), defined as the difference between the dual objective function value at time $t$ minus the primal objective function value at time $t$. This allows us to find close-to-optimal solutions for the consumers' price-response model~\eqref{eq:orig}. Thus, in the case when the duality gap is equal to zero, the measured consumption, if feasible, would be optimal in~\eqref{eq:orig}. In the case when the duality gap is greater than zero, the measured consumption would not be optimal. Intuitively, we attempt to find values for the parameters defining the block-wise utility function such that the measured consumption is as optimal as possible for problem~\eqref{eq:orig}.

Hence, the utility-refining problem consists in minimizing the sum of weighted duality gaps
\begin{subequations}\label{eq:primdual}
{\allowdisplaybreaks
\begin{IEEEeqnarray}{l}
\underset{ \begin{subarray}{c} a_{b,t},\lambda_{t}^{u}, \lambda_{t}^{d}, \\ \psi^{\overline{P}}_{t},\psi^{\underline{P}}_{t}, \underline{\psi}_{b,t}, \overline{\psi}_{b,t},\epsilon_t \end{subarray} }{\rm{Minimize}} \ \sum_{t \in \mathcal{T}} w_{t}\epsilon_t. \label{eq:obj_ref}
\end{IEEEeqnarray}
Note that we assign different weights to the duality gaps accrued in different time periods, in a way analogous to what we do with the absolute value of residuals in \eqref{eq:kktobj_lin1}. Objective function \eqref{eq:obj_ref} is subject to
\begin{IEEEeqnarray}{lLr}
&   \sum_{b \in \mathcal{B}} a_{b,1}x_{b,1}^{m'} - p_{1}\sum_{b \in \mathcal{B}} x_{b,1}  + \epsilon_1 =  \sum_{ b \in \mathcal{B}}  \left( \frac{\overline{P}_{1} - \underline{P}_{1}}{B} \right)  \overline{\psi}_{b,1}       &   \label{eq:objectivesprimdual} \IEEEeqnarraynumspace
\\
&   \sum_{b \in \mathcal{B}} a_{b,t}x_{b,t}^{m'} - p_{t}\sum_{b \in \mathcal{B}} x_{b,t}  + \epsilon_t = \sum_{ b \in \mathcal{B}}  \left( \frac{\overline{P}_{t} - \underline{P}_{t}}{B} \right)  \overline{\psi}_{b,t}     + & \nonumber
\\
&
 \left( r^{u}_{t} - \underline{P}_{t} + \underline{P}_{t-1} \right) \lambda_{t}^{u} +
\left( r^{d}_{t} + \underline{P}_{t} - \underline{P}_{t-1} \right)  \lambda_{t}^{d}
 \ \ t \in \mathcal{T}_ {-1}  \label{eq:objectivesprimdual2}
 &
\\
&\eqref{eq:stat1}-\eqref{eq:stat3}  \label{eq:stat_refine}&
\\
&a_{b,t} \geq a_{b+1,t}
\qquad\qquad \qquad\qquad t \in \mathcal{T} \label{eq:monotref}
\\
&\lambda_{t}^{u}, \lambda_{t}^{d}   \geq 0
\qquad\qquad\qquad\qquad \ \  t \in \mathcal{T}_{-1} \label{eq:feasref1}
\\
&\psi^{\overline{P}}_{t},\psi^{\underline{P}}_{t},\underline{\psi}_{b,t}, \overline{\psi}_{b,t}   \geq 0
\qquad\qquad  t \in \mathcal{T} \label{eq:feasref2}
\end{IEEEeqnarray}
}

\end{subequations}


The set of constraints \eqref{eq:objectivesprimdual2} constitutes the relaxed strong duality conditions, which express that the objective function of the original problem at time $t$, previously formulated in Equation \eqref{eq:orig}, plus the duality gap at time $t$, denoted by $\epsilon_t$, must be equal to the objective function of its dual problem also at time $t$. Equation \eqref{eq:objectivesprimdual} works similarly, but for $t=1$. The constraints relative to the dual of the original problem are grouped in \eqref{eq:stat_refine}.
Constraint \eqref{eq:monotref} requires that the estimated utility be monotonically decreasing. Finally, constraints \eqref{eq:feasref1} and \eqref{eq:feasref2} impose the non-negative character of dual variables.

%
\section{Case Study} \label{sec:studycase}
%
%

The proposed methodology to estimate the market bid that best captures the price-response of a pool of flexible consumers is tested using data from a real-life case study. The data relates to the Olympic Peninsula experiment, which took place in Washington and Oregon states between May 2006 and March 2007 \cite{Olympic}. The electricity price was sent out every fifteen minutes to 27 households that participated in the experiment. The price-sensitive controllers and thermostats installed in each house decided when to turn on and off the appliances, based on the price and on the house owner's preferences. 

\begin{figure}[h]
\centering
\includegraphics{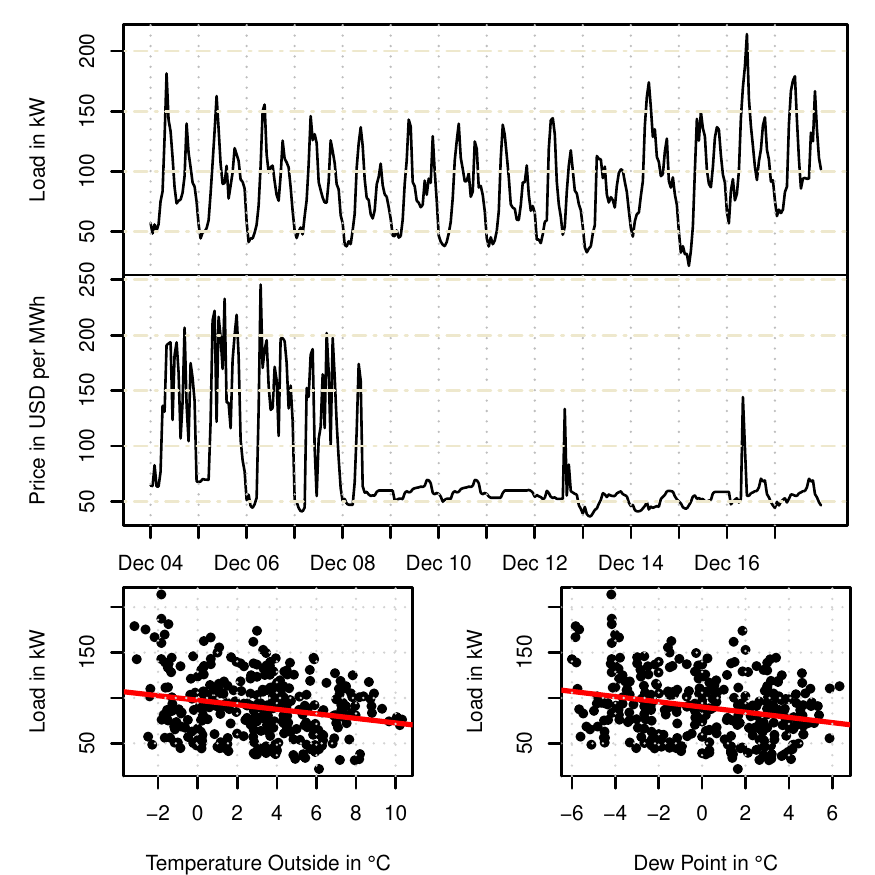}
\caption{ The upper and the middle plot show the load and the price, respectively. The bottom plots represent the load in the vertical axis versus the outside temperature and the dew point, on the left and on the right, respectively. The data shown spans from the 4th to the 18th of December.}
\label{fig:pricetempSep}
\end{figure}

For the case study, we use hourly measurements of load consumption, broadcast price, and observed weather variables, specifically, outside temperature, solar irradiance, wind speed, humidity, dew point and wind direction. Moreover, we include 0/1 feature variables to indicate the hour of the day, with one binary variable per hour (from 0 to 23), and one per day of the week (from 0 to 6). A sample of the dataset is shown in Figure \ref{fig:pricetempSep}, where the load is plotted in the upper plot, the price in the middle plot, and the load versus the outside temperature and the dew point in the bottom plots. The lines depicted in the bottom plots represent the linear relationship between the pairs of variables, and these are negative in both cases. The high variability in the price is also noteworthy: from the 1st to the 8th of December, the standard deviation of the price is 5.6 times higher than during the rest of the month (\$$67.9$/MWh versus \$$12.03$/MWh).


\subsection{Benchmark Models}


To test the quality of the market bid estimated by the proposed methodology, we quantify and assess the extent to which such a bid is able to \emph{predict} the consumption of the cluster of price-responsive loads. For the evaluation, we compare two versions of the inverse optimization scheme proposed in this paper with the Auto-Regressive model with eXogenous inputs (ARX) described in \cite{Corradi2012}. Note that this time series model was also applied by \cite{Corradi2012} to the same data set of the Olympic Peninsula project. All in all, we benchmark three different models:




\begin{description}
\item[\textbf{ARX},] which stands for Auto-Regressive model with eXogenous inputs \cite{MadsenBook}. This is the type of prediction model used in \cite{Dorini2012} and \cite{Corradi2012}. The consumption $x_{t}$ is modeled as a linear combination of past values of consumption up to lag $n$, $\boldsymbol X_{t-n} = \lbrace x_{t}, \ldots, x_{t-n} \rbrace$, and other explanatory variables $\boldsymbol Z_t = \lbrace Z_{t}, \ldots, Z_{t-n} \rbrace$. In mathematical terms, an ARX model can be expressed as $ x_{t} = \boldsymbol \vartheta_{x} \boldsymbol X_{t-n} + \boldsymbol\vartheta_{z} \boldsymbol Z_{t} + \epsilon_{t},$ with $\epsilon_{t} \sim \rm{N(0,}\sigma^2 \rm{)}$ and $\sigma^2$ is the variance.

\item[\textbf{Simple  Inv}]  This benchmark model consists in the utility-refining problem presented in Section \ref{sec:refine}, where the parameters of maximum pick-up and drop-off rates and consumption limits are computed from past observed values of consumption in a simple manner: we set the maximum pick-up and drop-off rates to the maximum values taken on by these parameters during the last seven days of observed data. All the features are used to explain the variability in the block-wise marginal utility function of the pool of price-responsive consumers: outside temperature, solar radiation, wind speed, humidity, dew point, pressure, and hour and week-day indicators. For this model, we use B=12 blocks of utility. This benchmark is inspired from the more simplified inverse optimization scheme presented in \cite{Keshavarz_imputinga} and \cite{chan2014generalized} (note, however, that neither \cite{Keshavarz_imputinga}, nor \cite{chan2014generalized} consider the possibility of leveraging auxiliary information, i.e., features, to better explain the data, unlike we do for the problem at hand).

\item[\textbf{Inv}]  This corresponds to the inverse optimization scheme with features that we propose, which runs following the two-step estimation procedure described in Section~\ref{sec:solution} with B=12 blocks of utility. Here we only use the outside temperature and hourly indicator variables as features. We re-parametrize weights $w_t$ with respect to a single parameter, called forgetting factor, and denoted as $E\geq0$, in the following manner:  $w_t = gap_t\left(\frac{t}{T}\right)^E$ for $t \in \mathcal{T}$ and $T$ being the total number of periods. The variable $gap$ indicates whether the observation was correctly measured ($gap = 1$) or not ($gap=0$). Parameter $E$ indicates how rapidly the weight drops (how rapidly the model forgets). When $E=0$, the weight of the observations is either 1 or 0 depending on the variable $gap$. As $E$ increases, the recent observations weight comparatively more than the old ones.



\end{description}

\subsection{Validation of the Model and Performance in December} \label{sec:validation}

In this subsection we validate the benchmarked models and assess their performance during the test month of December 2006.


For the sake of simplicity, we assume the price and the features to be known for the past and also for the future.
It is worth noticing, though, that the proposed methodology need not a prediction of the electricity price when used for bidding in the market and not for predicting the aggregated consumption of a cluster of loads. This is so because the market bid expresses the desired consumption of the pool of loads for any price that clears the market. The same cannot be said, however, for prediction models of the type of ARX, which would need to be used in combination with extra tools, no matter how simple they could be, for predicting the electricity price and for optimizing under uncertainty in order to generate a market bid.

\begin{figure}[h]
\centering
\includegraphics{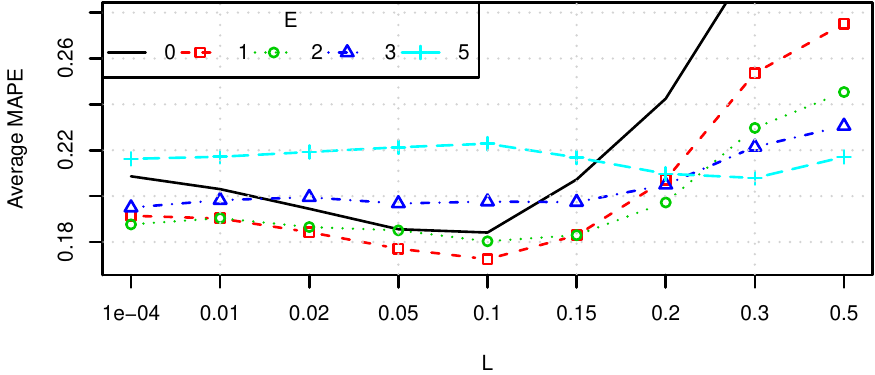}
\caption{Results from the validation of the input parameters $L$ and $E$, to be used during December.}
\label{fig:interaction}
\end{figure}

There are two parameters that need to be chosen before testing the models: the penalty parameter $L$ and the forgetting factor $E$. We seek a combination of parameters such that the prediction error is minimized. We achieve this by validating the models with past data, in a rolling-horizon manner, and with different combinations of the parameters $L$ and $E$.
The results are shown in Figure \ref{fig:interaction}. The MAPE is shown on the y-axis against the penalty $L$ in the x-axis, with the different lines corresponding to different values of the forgetting factor $E$. From this plot, it can be seen that a forgetting factor of  $E=1$ or $E=2$ yields a better performance than when there is no forgetting factor at all ($E=0$), or when this is too high ($E \geq3$). We conclude that selecting $L=0.1$ and $E=1$ results in the best performance of the model, in terms of the MAPE.

Once the different models have been validated, we proceed to test them. For this purpose, we first set the cross-validated input parameters to $L=0.1$ and $E=1$, and then, predict the load for the next day of operation in a rolling-horizon manner. In order to mimic a real-life usage of these models, we estimate the parameters of the bid on every day of the test period at 12:00 using historical values from three months in the past. Then, as if the market were cleared, we input the price of the day-ahead market (13 to 36 hours ahead) in the consumers' price-response model, obtaining a forecast of the consumption. Finally, we compare the predicted versus the actual realized consumption and move the rolling-horizon window to the next day repeating the process for the rest of the test period. Similarly, the parameters of the ARX model are re-estimated every day at 12:00, and predictions are made for 13 to 36 hours ahead.

\begin{figure}[h]
\centering
\includegraphics{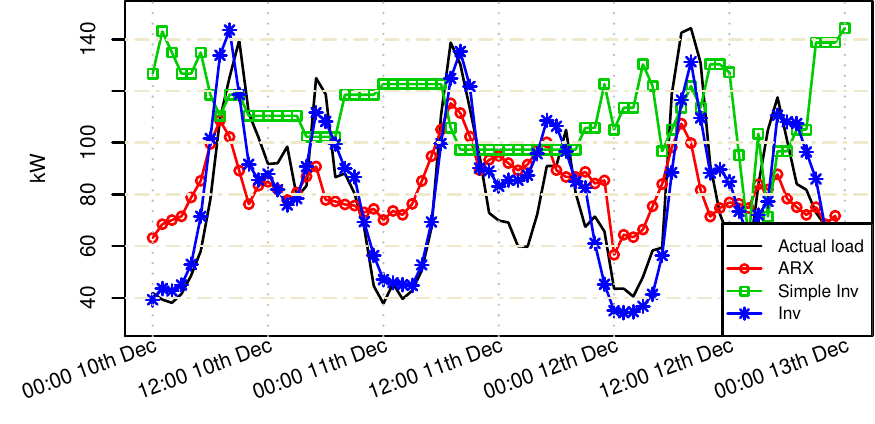}
\caption{Load forecasts issued by the benchmark models, and actual load, for the period between the 10th and the 13th of December.}
\label{fig:predicted}
\end{figure}

Results for a sample of consecutive days, from the 10th to the 13th of December, are shown in Figure \ref{fig:predicted}. The actual load is displayed in a continuous solid line, while the load predictions from the various benchmarked models are shown with different types of markers. First, note that the \textit{Simple Inv} model is clearly under-performing compared to the other methodologies, in terms of prediction accuracy. Recall that, in this model, the maximum and minimum load consumptions, together with the maximum pick-up and drop-off rates, are estimated from historical values and assumed to remain constant along the day, independently of the external variables (the features). This basically leaves the utility alone to model the price-response of the pool of houses, which, judging from the results, is not enough. The ARX model is able to follow the load pattern to a certain extent. Nevertheless, it is not able to capture the sudden decreases in the load during the night time or during the peak hours in the morning.
The proposed model (\textit{Inv}) features a considerably much better performance. It is able to follow the consumption pattern with good accuracy.

The performance of each of the benchmarked models during the whole month of December is summarized in Table \ref{tab:errors}. The first column shows the Mean Absolute Error (MAE), the second column provides the Root Mean Square Error (RMSE), and the third column collects the Mean Absolute Percentage Error (MAPE). The three performance metrics lead to the same conclusions: that the price-response models we propose, i.e., \textit{Inv}, perform better than the ARX model and the \textit{Simple Inv} model.
The results collated in Table \ref{tab:errors} also yield an interesting conclusion: that the electricity price is not the main driver of the consumption of the pool of houses and, therefore, is not explanatory enough to predict the latter. We conclude this after seeing the performance of the \textit{Simp Inv}, which is not able to follow the load just by modeling the price-consumption relationship by means of an utility function. The performance is remarkably enhanced when proper estimations of the maximum pick-up and drop-off rates and the consumptions bounds are employed.

\begin{table}
\centering
\caption{Performance measures for the three benchmarked models during December.}
 \begin{tabular}{cccc}
 \midrule
 & MAE & RMSE & MAPE   \\
  \cline{2-4}
    ARX    & 22.176 &  27.501 & 0.2750 \\
  	\textit{Simple Inv} 		& 44.437  & 54.576 & 0.58581 \\
    \textit{Inv} 		& 17.318  & 23.026 & 0.1899 \\
\bottomrule
\end{tabular}
\label{tab:errors}
\end{table}

%
%

\begin{figure}
\centering
\includegraphics{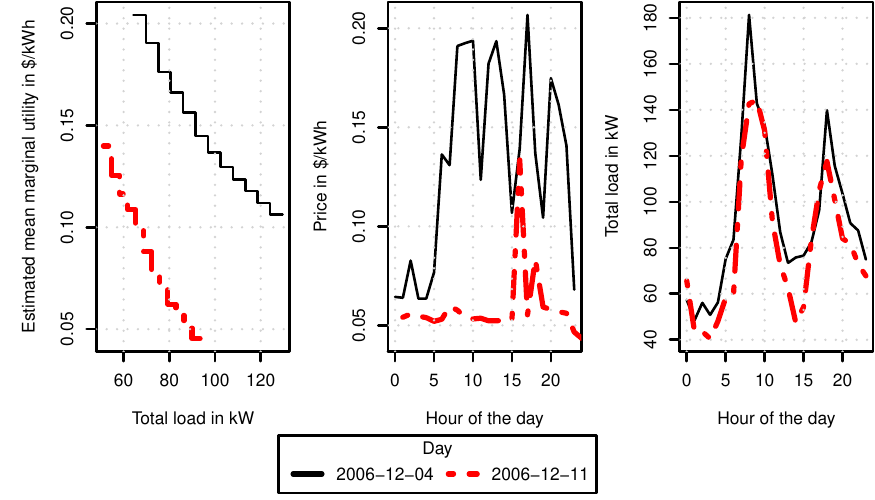}
\caption{Averaged estimated block-wise marginal utility function for the \textit{Inv} model (left panel), price in \$/kWh (middle panel), and load in kW (right panel). The solid lines represent data relative to the 4th of December. Dashed lines represent data relative to the 11th of December.}
\label{fig:util}
\end{figure}

The estimated block-wise marginal utility function, averaged for the 24 hours of the day, is shown in the left plot of Figure \ref{fig:util} for the \textit{Inv} model. The solid line corresponds to the 4th of December, when the price was relatively high (middle plot), as was the aggregated consumption of the pool of houses (right plot). The dashed line corresponds to the 11th of December and shows that the estimated marginal utility is lower, as is the price on that day.

\subsection{Performance During September and March}

In this section, we summarize the performance of the benchmarked models during September 2006 and March 2007.

In Table \ref{tab:errors2}, summary statistics for the predictions are provided for September (left side) and March (right side). The conclusions remain similar as the ones drawn for the month of December. The \textit{Inv} methodology consistently achieves the best performance during these two months as well.

\begin{table}[h]
\centering
\caption{Performance measures for the three benchmarked models}
 \begin{tabular}{cccc|ccc}
 \toprule
 & \multicolumn{3}{c}{September} &  \multicolumn{3}{c}{March} \\
 & MAE & RMSE & MAPE & MAE & RMSE & MAPE   \\
  \midrule
    ARX  &  7.649 & 9.829  & 0.2350 	   &		17.439 & 23.395 &  0.2509 \\
  	\textit{Simple Inv} 	& 14.263 & 17.8 &  0.4945 &		44.687 & 54.616 & 0.8365 \\
    \textit{Inv} 	& 5.719  & 8.582  & 0.1462 &		12.652  & 16.776 & 0.1952\\
    \bottomrule
\end{tabular}
\label{tab:errors2}
\end{table}

By means of cross-validation~\cite{james2013introduction}, we find that the user-tuned parameters yielding the best performance vary over the year. For September, the best combination is $L=0.3$, $E=0$, while for March it is $L=0.3$, $E=1$.

The optimized penalization parameter $L$ turns out to be higher in September and March than in December. This penalization parameter is highly related to the actual flexibility featured by the pool of houses. Indeed, for a high enough value of the penalty (say $L\geq0.4$ for this case study),  violating the complementarity conditions associated with the consumers' price-response model~\eqref{eq:orig} is relatively highly penalized. Hence, at the optimum, the slacks of the complementarity constraints in the relaxed estimation problem~\eqref{eq:penalization} will be zero or close to zero. When this happens, it holds at the optimum that $r^{u}_{t} = -r^{d}_{t}$ and $ \underline{P}_{t} = \overline{P}_{t}$. The resulting model is, therefore, equivalent to a linear model of the features, fit by least weighted absolute errors. When the best performance is obtained for a high value of $L$, it means that the pool of houses does not respond so much to changes in the price. On the other hand, as the best value for the penalization parameter $L$ decreases towards zero, the pool becomes more price-responsive: the maximum pick-up and drop-off rates and the consumption limits leave more room for the aggregated load to change depending on the price.

Because the penalization parameter is the lowest during December, we conclude that more flexibility is observed during this month than during September or March. The reason could be that December is the coldest of the months studied, with a recorded temperature that is on average 9.4$^{\circ}$C lower, and it is at times of cold whether when the electric water heater is used the most.

\section{Summary and Conclusions}

We consider the market-bidding problem of a pool of price-responsive consumers. These consumers are, therefore, able to react to the electricity price, e.g., by shifting their consumption from high-price hours to lower-price hours. The total amount of electricity consumed by the aggregation has to be purchased in the electricity market, for which the aggregator or the retailer is required to place a bid into such a market. Traditionally, this bid would simply be a forecast of the load, since the load has commonly behaved inelastically. However, in this paper, we propose to capture the price-response of the pool of flexible loads through a more complex, but still quite common market bid that consists of a stepwise marginal utility function, maximum load pick-up and drop-off rates, and maximum and minimum power consumption, in a manner analogous to the energy offers made by power producers.

We propose an original approach to estimate the parameters of the bid based on inverse optimization and bi-level programming. Furthermore, we use auxiliary variables to better explain the parameters of the bid. The resulting non-linear problem is relaxed to a linear one, the solution of which depends on a penalization parameter. This parameter is chosen by cross-validation, proving to be adequate from a practical point of view.


For the case study, we used data from the Olympic Peninsula project to asses the performance of the proposed methodology. We have shown that the estimated bid successfully models the price-response of the pool of houses, in such a way that the mean absolute percentage error incurred when using the estimated market bid for predicting the consumption of the pool of houses is kept in between 14\% and 22\% for all the months of the test period.


We envision two possible avenues for improving the proposed methodology. The first one is to better exploit the information contained in a large dataset by allowing for non-linear dependencies between the market-bid parameters and the features. This could be achieved, for example, by the use of B-splines. The second one has to do with the development of efficient solution algorithms capable of solving the exact estimation problem within a reasonable amount of time, instead of the relaxed one. This could potentially be accomplished by decomposition and parallel computation.

%
%

\ifCLASSOPTIONcaptionsoff
  \newpage
\fi



%

\bibliographystyle{IEEEtran}
%

\appendix

Next we show how to formulate robust constraints to ensure that the estimated minimum consumption be always equal to or lower than the estimated maximum consumption. At all times, and for all plausible realizations of the external variables, we want to make sure that:
\begin{align}
\underline{P}^0 + \sum_{i \in \mathcal{I}} \alpha_{i}^{\underline{P}} Z_{i,t} \leq \overline{P}^0  + \sum_{i \in \mathcal{I}}  \alpha_{i}^{\overline{P}} Z_{i,t}, &  \ \ \ t \in \mathcal{T}, \ \forall Z_{i,t}. \label{eq:robustbound1}
\end{align}

If \eqref{eq:robustbound1} is not fulfilled, problem~\eqref{eq:orig} is infeasible (and the market bid does not make sense). Assuming we know the range of possible values of the features, i.e.,  $Z_{i,t} \in [\overline{Z}_i , \underline{Z}_i]$, \eqref{eq:robustbound1} can be rewritten as:

\begin{align}
\underline{P}^0 -  \overline{P}^0 +    \underset{ \begin{subarray}{c} Z_{i,t}' \\ \rm{s.t.} \ \underline{Z}_i \leq Z_{i,t}' \leq \overline{Z}_{i} \\ {i \in \mathcal{I}} \end{subarray} }{ \rm{Maximize}} \left\lbrace
         \sum_{i \in \mathcal{I}} ( \alpha_{i}^{\underline{P}} - \alpha_{i}^{\overline{P}}) Z_{i,t}'
      \right \rbrace \leq 0, \ t \in \mathcal{T}. \label{eq:robustbound}
\end{align}

Denote the dual variables of the upper and lower bounds of $Z_{i,t}'$ by $\overline{\phi}_{i,t}$ and $\underline {\phi}_{i,t}$ respectively. Then, the robust constraint~\eqref{eq:robustbound} can be equivalently reformulated as:

\begin{subequations}
\begin{align}
&\overline{P}^0 - \underline{P}^0 +\sum_{i \in \mathcal{I}} (\overline{\phi}_{i,t} \overline{Z}_i -  \underline{\phi}_{i,t} \underline{Z}_i) \leq 0 & t \in \mathcal{T}\label{eq:robust1}
\\
&\overline{\phi}_{i,t} - \underline{\phi}_{i,t} = \alpha_{i}^{\overline{P}} - \alpha_{i}^{\underline{P}}  & i \in \mathcal{I} ,t \in \mathcal{T}
\\
&\overline{\phi}_{i,t},  \underline{\phi}_{i,t} \geq 0 & i \in \mathcal{I},t \in \mathcal{T}.\label{eq:robust3}
\end{align}
\end{subequations}

Following the same reasoning, one can obtain the set of constraints that guarantees the non-negativity of the lower bound and consistent maximum pick-up and drop-off rates. We leave the exposition and explanation of these constraints out of this paper for brevity.

\end{document}